%% file: SC_ddfv_vs_hfv.tex
\documentclass{svproc}
\usepackage{url}

\usepackage{amsfonts}
\usepackage{amssymb}
\usepackage{amsmath}
\usepackage{amssymb}
\usepackage{mathrsfs}
\usepackage{graphicx}
\usepackage{caption} 
\usepackage{wrapfig}
\usepackage[dvipsnames]{xcolor}
\usepackage{dsfont}
\usepackage{empheq}
\usepackage{thmtools}
\usepackage{thm-restate}
\usepackage[normalem]{ulem}
\usepackage{subcaption}
\usepackage[mathscr]{eucal}

\usepackage{hyperref}
\hypersetup{
  unicode = true,                           
  pdftoolbar = true,                        
  pdfmenubar = true,                        
  pdffitwindow = true,                      
  pdfauthor = {S. Krell and J. Moatti},              
  pdftitle = {Structure-preserving schemes for drift-diffusion systems on general meshes: DDFV vs HFV},         
  pdfnewwindow = true,                      
  colorlinks = true,                        
  linkcolor = blue,                         
  citecolor = magenta,                      
  filecolor = violet,                       
  urlcolor = violet                         
}

\usepackage{mathptmx}       
\usepackage{helvet}         
\usepackage{courier}        
\usepackage{type1cm}        
\usepackage{makeidx}         
\usepackage{graphicx}        
\usepackage{multicol}        
\usepackage[bottom]{footmisc}
\usepackage{amsmath, amsfonts, amssymb}
\usepackage{color}
\usepackage{tikz}

\usepackage{graphicx}
\usepackage{tikzscale}
\usepackage{paralist}
\usepackage{pgfplots}
\pgfplotsset{every axis/.append style={
                    label style={font=\Large},
                    tick label style={font=\Large},
                    legend style={font=\Large}
                    }}

\usepackage{dsfont}

\usepackage{stmaryrd}
\usepackage{wrapfig}

\DeclareMathOperator{\e}{e}
\DeclareMathOperator{\divergence}{div}

\newcommand{\Disc}{ \EuScript{D}}

\newcommand{\E}{\mathcal{E}}
\newcommand{\s}{\sigma}
\newcommand{\V}{\underline{V}}

\renewcommand{\u}{\underline{u}}
\renewcommand{\v}{\underline{v}}
\newcommand{\w}{\underline{w}}

\newcommand{\one}{\underline{1}}

\newcommand{\Nd}{\underline{N}}
\newcommand{\Pd}{\underline{P}}

\newcommand{\phid}{\underline{\phi}}

\usepackage[justification=centering]{caption}
\input{macroDDFV.tex}

\definecolor{MyGreen}{RGB}{54,165,54}
\newcommand{\JM}[1]{{#1}}
\definecolor{forestgreen}{rgb}{0.13,0.54,0.13}

\begin{document}
\mainmatter              

\title{Structure-preserving schemes for drift-diffusion systems on general meshes: DDFV vs HFV}
\titlerunning{Comparison of DDFV and HFV structure preserving schemes}
\author{Stella Krell\inst{1} \and Julien Moatti\inst{2}}
\authorrunning{Stella Krell and Julien Moatti} 
%
\tocauthor{Ivar Ekeland, Roger Temam, Jeffrey Dean, David Grove,
Craig Chambers, Kim B. Bruce, and Elisa Bertino}
\institute{	Universit\'e C\^ote d'Azur, CNRS, Inria, LJAD, France, \\ 
				\email{stella.krell@univ-cotedazur.fr}
	\and
			Inria, Univ. Lille, CNRS, UMR 8524 - Laboratoire Paul Painlev\'e, F-59000 Lille, France \\
				\email{julien.moatti@inria.fr}
} 

\maketitle              

\begin{abstract}
We made a comparison between a Discrete Duality Finite Volume (DDFV) scheme and a Hybrid Finite Volume (HFV) scheme for a drift-diffusion model with mixed boundary conditions on general meshes. 
Both schemes are based on a nonlinear discretisation of the convection-diffusion fluxes, which ensures the positivity of the discrete densities. 
We investigate the behaviours of the schemes on various numerical test cases.
\end{abstract}


\section{Motivation}

%
We are interested in the numerical discretization of drift-diffusion model. Let $\Omega$ be a polygonal connected open bounded subset of $\R^2$, whose boundary $\Gamma=\partial \Omega$ is divided into two parts $\Gamma=\Gamma^D\cup\Gamma^N$ with ${\rm m}(\Gamma^D)>0$. The problem writes:
\begin{equation} \label{pb:DD:evol}	
	\left \lbrace
	\begin{aligned}
        \partial _ t N - \divergence ( \nabla N - N \nabla \phi   ) 	& = 0 	&&\text{ in } \R_+ \times \Omega, \\
		\partial _ t P - \divergence (\nabla P  + P \nabla \phi  ) 	& = 0 		&&\text{ in } \R_+ \times \Omega, \\
		- \lambda^2 \divergence ( \nabla \phi  ) &=  C + P - N  &&\text{ in } \R_+ \times \Omega, \\
		N =  N^D, \ P =  P^D \text{ and } \phi &=  \phi^D &&\text{ on }  \R_+ \times \Gamma^D, \\
		( \nabla N - N \nabla \phi ) \cdot n = (  \nabla P  + P \nabla \phi  ) \cdot n 
		& =  \nabla \phi  \cdot n  = 0 &&\text{ on } \R_+ \times \Gamma^N,\\		
	 	N(0, \cdot)  =  N^{in} \text{ and } \ P(0, \cdot) &=  P^{in}   &&\text{ in } \Omega, 
    \end{aligned}
	\right.
\end{equation}
where $n$ denotes the unit normal vector to $\partial \Omega$ pointing outward $\Omega$.
Regarding the data, 
(i) the parameter $\lambda >0$ is the rescaled Debye length of the system, which accounts for the nondimensionalisation (relevant values of this parameter can be very small, inducing some stiff behaviours), 
(ii) the initial conditions $N^{in}$ and $P^{in}$ belong to $L^\infty(\Omega)$ and are positive, 
(iii) the doping profile $C$ is in $L^\infty(\Omega)$, and characterises the semiconductor device used. 
In the following, we also assume that the boundary conditions are the trace of some $H^1$ function on $\Omega$, such that the following relation holds:
\begin{equation}\label{eq:compcond}
	\log(N^D) - \phi^D = \alpha_N \text{ and } \log(P^D) + \phi^D = \alpha_P \text{ on }  \Gamma^D, 
\end{equation}
where $\alpha_N$ and $\alpha_P$  are two real constants. It follows that $N^D$ and $P^D$ are positive. 
%

%
The solution to \eqref{pb:DD:evol} enjoys some natural physical properties: the densities $N$ and $P$ are positive for all time, and the solution converges exponentionaly fast towards some thermal equilibrium $(N^e,P^e,\phi^e)$ -which is a stationary solution to \eqref{pb:DD:evol}- where
$N^e = \e^{\alpha_N + \phi^e}$, $P^e =\e^{\alpha_P - \phi^e}$ and $\phi^e$ is the solution to the Poisson-Boltzmann equation 
\begin{equation} \label{pb:PoissonBoltz}
	\left \lbrace
	\begin{aligned}
		- \lambda^2 \divergence ( \nabla \phi^e  ) & =  C +  \exp(\alpha_P - \phi^e ) - \exp(\alpha_N + \phi^e )  &&\text{ in } \Omega, \\
		\phi^e &=  \phi^D \text{ on }  \Gamma^D 
		 \qquad \text{ and } \qquad 	
		\nabla \phi^e  \cdot n  = 0 &&\text{ on }  \Gamma^N.
	\end{aligned}
	\right.
\end{equation}
Relation \eqref{eq:compcond} is a compatibility condition in order to ensure the existence of the thermal equilibrium \eqref{pb:PoissonBoltz}.
%
When designing numerical schemes for \eqref{pb:DD:evol}, it is crucial to ensure that the scheme preserves these properties at the discrete level.
This structure preserving feature is ensured by classical TPFA schemes on admissible orthogonal meshes (see \cite{BCCH:17}).
Unfortunately, these schemes cannot be used on general meshes. 
Following the ideas introduced in \cite{CaGui:17}, a nonlinear positivity preserving DDFV scheme for Fokker-Planck equations has been introduced in \cite{CCHK18}.
In the spirit of these works, a nonlinear structure preserving HFV scheme was introduced and partially analysed in \cite{Moa:22}. 
%
The aim of this paper is to introduce a nonlinear structure preserving DDFV scheme for \eqref{pb:DD:evol} based on the scheme of \cite{CCHK18} and to compare it numerically with the HFV scheme of \cite{Moa:22}. 
\section{Descriptions of the schemes}
The schemes used here are based on the same nonlinear strategy, introduced in \cite{CaGui:17}, consisting in the reformulation of the convection-diffusion fluxes: 
\[
	\nabla N - N \nabla \phi = N \nabla \left( \log(N)  - \phi \right ) 
	\text{ and } 
	\nabla P + P \nabla \phi = P \nabla \left( \log(P)  + \phi \right ) .
\]
At the discrete level, both schemes relie on discrete gradients operators to approximate the continuous gradients.
The major issue lies in the discretisation of the prefactors $P$ and $N$, which will be handled by local reconstruction operators.
\JM{
This discretisation strategy is a way of ensuring (at the theoretical level) the positivity of the discrete densities.
We refer to \cite[Theorem 1]{Moa:22} (HFV scheme for drift-diffusion system) and \cite[Theorem 2.1]{CCHK18} (DDFV scheme for a single advection-diffusion equation) for proofs of this statement. We also refer the reader to these proofs for more insight about the reconstruction operators.
}
Both schemes are based on a backward Euler discretisation in time. To fix ideas, we will use a constant time step $\Delta t > 0$.
For more precise descriptions and statements about the schemes and the meshes, we refer to \cite{CCHK18} (DDFV) and \cite{Moa:22} (HFV).

\begin{remark}[Generalisation to anisotropic models]
In this paper, we consider isotropic convection-diffusion equations for the charges carriers for the sake of brevity. 
One could add anisotropic diffusion tensors and consider the framework described in \cite{Moa:22}.  
\end{remark}

Both schemes rely on a spatial discretisation (or mesh) of the domain $\Omega$.
The (primal interior) mesh $\M$ is a partition of $\Omega$ in polygonal control volumes (or cells).
We let $\dr\M$ be the set of boundary edges, seen either as degenerate control volumes (DDFV framework) or as edges (HFV framework).
The primal mesh $\overline{\M}$ is defined as the reunion of $\M$ and $\dr\M$.
Given a cell $K\in \overline{\M}$, we fix a point $\xk \in K$, called the center of $K$. 
For all neighboring primal cells $\k$ and $\l$, we assume that $\dr\k\cap\dr\l$ is a segment, corresponding to an internal edge of the mesh $\M$, denoted by $\sigma=\k\vert\l$ and we let $\Ee_{int}$ be the set of such edges.
We denotes by $\E = \E_{int} \cup \dr\M$ the set of all (internal and exterior) edges of the mesh, and define $\E_K$ the set of edges of the cell $K \in \M$. For any $K \in \M$ and $\s \in \E_K$, we define $\nksig$ as the unit normal to $\sigma$ outward $\k$.
\newline
Given any measurable $X \subset \R^2$, we denote by $m_X$ the measure of the object $X$.
\subsection{The DDFV scheme}
%
In order to define the DDFV scheme, we need to introduce two other meshes: 
the dual mesh denoted $\overline{\Mie} $ and the diamond mesh denoted $\DD$ (see \cite{CCHK18} for more details).
The dual mesh $\overline{\Mie}$ is also composed of interior dual mesh $\Mie$ (corresponding of cells around vertex in $\Omega$) and of boundary dual mesh $\dr\Mie$ (corresponding of cells around vertex on $\dr\Omega$).
For any vertex $\xke$ of the primal mesh satisfying $\xke\in \Omega$, we define a polygonal control volume $\ke$ by connecting all the centers of the primal cells sharing $\xke$ as vertex.
For any vertex $\xke\in \partial \Omega$, we define a polygonal control volume $\ke$ by connecting the centers $\xk $ of the interior primal cells and the midpoints of the boundary edges sharing $\xke$ as vertex and $\xke$. 
We define the set $\Ee^*_{int}$ of internal edges of the dual mesh similarly as $\Ee_{int}$. We denote by $\nkesige$ the unit normal to $\sigma^*$ outward $\ke$.
For each couple $(\sigma,\sigma^*)\in\Ee_{}\times\Ee^*_{int}$ such that $\sigma=[\xke,\xle]$ and $\sigma^*=\ke\vert\le$, we define the quadrilateral diamond $\Dsig$ whose diagonals are $\sigma$ and $\sigma^*$ (if $\sigma\subset\partial\Omega$, it degenerates into a triangle).
The set of the diamonds defines the diamond mesh $\DD$ , which is a partition of $\Omega$.
\vspace{-0.2cm}
\begin{figure}[htb]
\begin{center}
\begin{tikzpicture}[scale=.95]
\node[rectangle,scale=0.8,fill=black!50] (xle) at (0,0) {};
\node[circle,draw,scale=0.5,fill=black!5] (xl) at (2,1.3) {};
\node[rectangle,scale=0.8,fill=black!50]  (xke) at (0,4) {};
\node[circle,draw,scale=0.5,fill=black!5] (xk) at (-2,2.3) {};
\draw[line width=1pt] (xle)--(xke);
\draw[dashed, line width=1pt] (xk)--(xl);
\draw[dash pattern=on 2pt off 3pt on 6pt off 3pt,line width=2pt] (xk)--(xke)--(xl)--(xle)--(xk);

\node[yshift=-8pt] at (xle){$\xle$};
\node[yshift=8pt] at (xke){$\xke$};
\node[xshift=10pt] at (xl){$\xl$};
\node[xshift=-10pt] at (xk){$\xk$};

\draw[->,line width=1pt] (-2+3.*0.4,2.3-3.*0.1)--(-2+4.8*0.4,2.3-4.8*0.1);
\draw[->,line width=1pt] (-2+3.*0.4,2.3-3.*0.1)--(-2+3.*0.4-1.8*0.1,2.3-3.*0.1-1.8*0.4);
\draw[->,line width=1pt] (0,2.7)--(0,2.);
\draw[->,line width=1pt] (0,2.7)--(.7,2.7);
\node[right] at (0,2.3){$\tkele$};
\node[right] at (0.1,2.95){$\nksig$};
\node at (-0.3,1.6){$\tkl$};
\node at (-0.5,1.2){$\nkesige$};

\draw[line width=1pt]  (2.2,4)--(2.8,4);
\node[right,xshift=8] at (2.8,4){$\sigma=\k\vert\l$, edge of the primal mesh};
\draw[dashed, line width=1pt] (2.2,3.5)--(2.8,3.5);
\node[right,xshift=8] at (2.8,3.5){$\sigma^*=\ke\vert\le$, edge of the dual mesh};
\draw[dash pattern=on 2pt off 3pt on 6pt off 3pt,line width=2pt] (2.2,3)--(2.8,3);
\node[right,xshift=8] at (2.8,3.){Diamond $\Dsig$};
\node[rectangle,scale=0.8,fill=black!50] at (2.8,2.5) {};
\node[circle,draw,scale=0.5,fill=black!5] at (2.8,2.) {};
\node[right,xshift=8] at (2.8,2.5){Vertices of the primal mesh};
\node[right,xshift=8] at (2.8,2.){Centers of the primal mesh};
\node [right]at (0,1.5) {$x_\Ds$};

\node[rectangle,scale=0.8,fill=black!50] (xle2) at (9,0) {};
\node[circle,draw,scale=0.5,fill=black!5] (xl2) at (9,2) {};
\node[rectangle,scale=0.8,fill=black!50]  (xke2) at (9,4) {};
\node[circle,draw,scale=0.5,fill=black!5] (xk2) at (7.5,2.3) {};
\draw[line width=1pt] (xle2)--(xke2);
\draw[dashed, line width=1pt] (xk2)--(xl2);
\draw[dash pattern=on 2pt off 3pt on 6pt off 3pt,line width=2pt] (xle2)--(xk2)--(xke2);

\node[yshift=-8] at (xle2){$\xle$};
\node[yshift=8] at (xke2){$\xke$};
\node[xshift=10] at (xl2){$\xl$};
\node[xshift=-10] at (xk2){$\xk$};

\end{tikzpicture}
\end{center}
\vspace{-0.5cm}
\caption{Definition of the diamonds $\Dsig$ and related notations.}\label{fig_diamonds}
\end{figure}
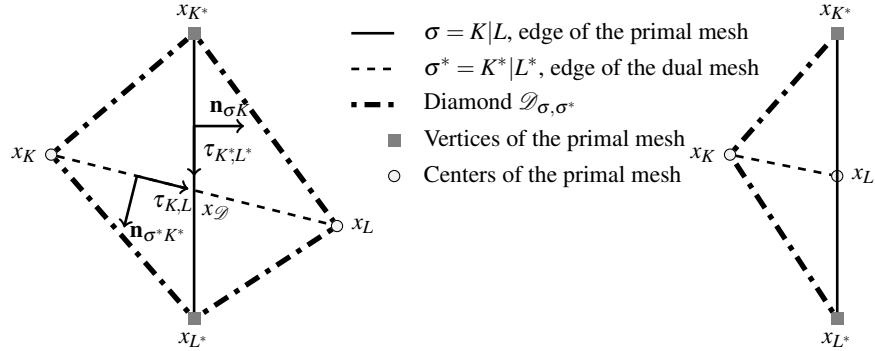
Finally, the DDFV mesh is made of  $\T=(\overline{\M},\overline{\Mie})$ and $\DD$.
%
%

We now introduce the space of scalar fields which are associated to each primal and dual cell $\Rt$, and space of vector fields constant on the diamonds $\RdDD$:
\vspace{-0.2cm}
\[
 \ut\in\Rt\Longleftrightarrow  \ut =\left(\left(\uk\right)_{K\in{\overline \M}},\left(\uke\right)_{K^\ast\in{\overline \Mie}}\right)
 	\text{ and } 
  \xib_\DD \in \RdDD\Longleftrightarrow \xib_\DD =\left(\xib_\Ds\right)_{\Ds\in\DD}.
\]
To enforce Dirichlet boundary conditions, we introduce the set of Dirichlet boundary primal and dual cells:
$
\partial \mathfrak{M}_D = \{ \k \in \partial \mathfrak{M}: \k \subset  \Gamma_D\}
	$ and  $
\partial \mathfrak{M}^*_D = \{ \ke \in \partial \mathfrak{M^*}: \xke \in  {\overline \Gamma_D}\},
$
and, for a given $v\in C(\Gamma^D)$, we define 
\begin{equation*}
\textsc{E}_{v}^{\Gamma_D}= \lbrace \ut \in \Rt \mid
\forall \k \in \partial \mathfrak{M}_{D}, \, \uk=v(\xk) \text{ and } 
\forall \ke \in \partial \mathfrak{M}^*_{D}, \, \uke=v(\xke) \rbrace.
\end{equation*}
We also define discrete bilinear forms on $\Rt$ and $\RdDD$ by
\vspace*{-0.2cm}
\begin{align*}
&\dsp\left\llbracket\vt,\ut\right\rrbracket_\T&&=&&\frac{1}{2}\left(\sumpri\mk\uk \vk+\sumdua \mke\uke\vke
\right), \ \ \ \forall (\ut,\vt) \in \left (\Rt \right) ^2,\\
&\dsp\left(\xib_\DD,\phib_\DD\right)_{\DD}&&=&&\sumdiam\md\ 
 \xib_\Ds\cdot \phib_\Ds,\ \ \ \forall (\xib_\DD,\phib_\DD) \in \left ( \RdDD \right )^2.
\end{align*}
%
The DDFV method is based on the definition of a discrete gradient operator $\nabla^\DD : \Rt \to \RdDD$, defined by 
$\gradDD \ut =\dsp\left(\gradD \ut\right)_{\Ds\in\DD}$, where
\vspace{-0.3cm}
\begin{equation}
	\gradD \ut =\frac{1}{2 m_\Ds} \left(  \msig (\ul-\uk)\nksig + \msige(\ule-\uke)\nkesige \right) \quad \forall  \Ds\in\DD.
\end{equation}
Finally, we introduce a reconstruction operator on diamonds $r^\DD$. It is a mapping from  $\R^\T$ to  $\R^\DD$ defined  for all  $\ut\in \Rt$ by  $r^\DD \ut =\left(r^\Ds \ut\right)_{\Ds\in\DD}$, 
where for $\Ds\in\DD$, whose vertices are $x_K$, $x_L$, $x_\ke$, $x_\le$,
$
r^\Ds \ut=\frac{1}{4}(\uk+\ul+\uke+\ule).
$
One can now introduce a DDFV discretisation of $(u,w,v) \mapsto \int_\Omega u \nabla w \cdot \nabla v$, defined by 
\[
	T_{\DD} : (\ut, \wt,  \vt ) \mapsto
		\sumdiam \md r^\Ds \ut \; \gradD\wt \cdot  \gradD\vt.\vspace{-0.3cm}
\]
Now, we first discretise the data by taking the mean values of $N^{in}$, $P^{in}$ and $C$  on the primal and dual cells, which define  $\Nt^0$, $\Pt^0$ and $C_\mathcal{T}$.
Then, for all $n\geq 0$, we look for $(\Nt^{n+1},\Pt^{n+1},\pott^{n+1})\in \textsc{E}_{N^D}^{\Gamma_D}\times\textsc{E}_{P^D}^{\Gamma_D}\times\textsc{E}_{\phi^D}^{\Gamma_D}$ solution to:
\begin{subequations}\label{sch:DDFV}
\begin{gather}
	\Bigl\llbracket\dsp\frac{\Nt^{n+1}-\Nt^n}{\Delta t}, \vt \Bigl\rrbracket_\T
			+T_{\DD}(\Nt^{n+1}, \log (\Nt^{n+1}) - \pott^{n+1}, \vt )=0 
			\quad \forall \vt \in\textsc{E}_0^{\Gamma_D},\label{sch:DDFV:N}\\
	\Bigl\llbracket\dsp\frac{\Pt^{n+1}-\Pt^n}{\Delta t}, \vt \Bigl\rrbracket_\T
			+T_{\DD}(\Pt^{n+1}, \log (\Pt^{n+1})+\pott^{n+1} ,\vt)=0 
			\quad \forall \vt \in\textsc{E}_0^{\Gamma_D},\label{sch:DDFV:P}\\
	\lambda^2 \left( \gradD  \pott^{n+1} , \gradD  \vt \right)_{\DD}
 			=\Bigl\llbracket C_\T+\Pt^{n+1}-\Nt^{n+1}, \vt \Bigl\rrbracket_\T 
 			\quad \forall \vt \in\textsc{E}_0^{\Gamma_D}.\label{sch:DDFV:phi}
\end{gather}
\end{subequations}
In \eqref{sch:DDFV:N} and \eqref{sch:DDFV:P}, we use the notation $ \log (\ut) = \left(\left( \log(u_K) \right)_{K\in{\overline \M}}, \left( \log(u_{K^\ast})\right)_{K^\ast\in{\overline \Mie}}\right)$.
\subsection{The HFV scheme} 
In order to define the HFV scheme, we need to introduce a pyramidal submesh. To do so, one has to assume that each cell $K \in \M$ is star-shaped with respect to its center $x_K$ (we recall that $x_K$ is not necessarily the barycentre of $K$). 
We then define $P_{K,\s}$ as the pyramid (triangle) of base $\s$ and apex $x_K$. 
Given any $\s\in \E$, we denote by $\overline{x}_\s$ the barycentre of $\s$, and by $d_{K,\s}$ the euclidean distance between $\s$ and $x_K$.
Finally, we define the hybrid discretisation (or mesh) as $\Disc = ( \M, \E )$.

%
We now introduce the space of discrete (scalar) hybrid unknowns $\V_\Disc$:
\vspace{-0.1cm}
\[
  \u_\Disc \in \V_\Disc \Longleftrightarrow  \u_\Disc 	
  	= \left( \left(\uk\right)_{K\in{\M}},\left(u_\s\right)_{\s \in \E}\right), 
\]
where the $u_K \in \R$ are the cell unknowns and the $u_\s \in \R$ are the edges unknowns (approximation of the trace of the solutions on the edges).
To enforce Dirichlet boundary conditions, for a given $v\in C(\Gamma^D)$, we define 
\vspace{-0.17cm}\begin{equation*}
\V_{\Disc,v}^{\Gamma_D}= \lbrace \u_\Disc \in \V_\Disc \mid
\forall \s \in \partial \mathfrak{M}_{D}, \, u_\s =v(\overline{x}_\s)  \rbrace.
\end{equation*}
As for the DDFV framework, we define a bilinear form on $\V_\Disc$, discrete counterpart of the inner product on $L^2(\Omega)$ as
\[
\dsp\left\llbracket\u_\Disc,\v_\Disc\right\rrbracket_\M = \sumpri\mk\uk \vk, \ \ \ \forall (\u_\Disc,\v_\Disc)\in  \V_\Disc^2.
\]
\vspace*{-0.2cm}

%
The HFV method is based on the definition of a discrete gradient operator $\nabla_\Disc : \V_\Disc \to (\R^2)^\Omega$ which maps discrete hybrid unknowns onto piecewise constant functions on the pyramidal submesh. More precisely, given $\v_\Disc \in \V_\Disc$, $K \in \M$ and $\s \in \E_K$, \vspace*{-0.12cm}
\[
  {\nabla_{\Disc}\v_{\Disc}}_{\mid P_{K,\s}} =G_K \v_\Disc + S_{K,\s} \v_\Disc,
\]
where, for some $\eta >0$, the consistent and stabilisation parts of the gradient are given by 
\vspace{-0.12cm}
\[
	G_K\v_\Disc=\frac{1}{\mk}\sum_{{\s'}\in\E_K}m_{\s'}v_{\s'}n_{K,\s'}
	\text{ and }	
	S_{K,\s}\v_\Disc=\frac{\eta}{d_{K,\s}}\big(v_\s-v_K-G_K\v_K\cdot(\overline{x}_\s-x_K)\big)n_{K,\s}.
\]\vspace{-0.12cm}
One can now define the discrete counterpart of $(u,v) \mapsto \int_\Omega \nabla u \cdot \nabla v$ as 
\[
	a_\Disc : (\u_\Disc, \v_\Disc) \mapsto \int_\Omega \nabla_{\Disc}\u_{\Disc} \cdot \nabla_{\Disc}\v_{\Disc}.
\]
Finally, we introduce as previously local reconstruction operators on cells $r^K : \V_\Disc \to \R $, such that for any $u_\Disc \in \V_\Disc$,
$ \displaystyle 
	r^K( \u_\Disc)  = \frac{1}{|\E_K|}\sum_{\s \in \E_K} \frac{u_K + u_\s}{2}, 
$
where $|\E_K|$ is the cardinal of the finite set $\E_K$.
One can now introduce a HFV discretisation of $(u,w,v) \mapsto \int_\Omega u \nabla w \cdot \nabla v$, defined by 
\vspace{-0.3cm}
\[
	T_{\Disc} : (\u_\Disc, \w_\Disc,  \v_\Disc ) \mapsto
		\sum_{K \in \M} r^K( \u_\Disc) \; \int_K \grad_\Disc \w_\Disc \cdot  \grad_\Disc \v_\Disc.
\]
%
We now discretise the data by taking the mean values of $N^{in}$, $P^{in}$ and $C$ on the cells and edges, which define  $\Pd_\Disc^0$, $\Nd_\Disc^0$ and $\underline{C}_\Disc$. 
Then, for all $n\geq 0$, we look for 
$	(\Nd_\Disc^{n+1},\Pd_\Disc^{n+1},\phid_\Disc^{n+1})
		\in  \V_{\Disc,N^D}^{\Gamma_D} \times \V_{\Disc,P^D}^{\Gamma_D}  \times \V_{\Disc,\phi^D}^{\Gamma_D}$ solution to:
\begin{subequations}\label{sch:HFV}
\begin{gather}
	\Bigl\llbracket\dsp\frac{\Nd_\Disc^{n+1}-\Nd_\Disc^{n}}{\Delta t}, \v_\Disc \Bigl\rrbracket_\M
			+T_{\Disc}(\Nd_\Disc^{n+1}, \log (\Nd_\Disc^{n+1}) - \phid_\Disc^{n+1},\v_\Disc  )=0 
			\quad \forall \v_\Disc\in \V_{\Disc,0}^{\Gamma_D},\label{sch:HFV:N}\\
	\Bigl\llbracket\dsp\frac{\Pd_\Disc^{n+1}-\Pd_\Disc^{n}}{\Delta t}, \v_\Disc \Bigl\rrbracket_\M
			+T_{\Disc}(\Pd_\Disc^{n+1}, \log (\Pd_\Disc^{n+1}) + \phid_\Disc^{n+1},\v_\Disc  )=0 
			\quad \forall \v_\Disc \in \V_{\Disc,0}^{\Gamma_D},\label{sch:HFV:P}\\
	\lambda^2 a_\Disc \left( \phid_\Disc^{n+1} , \v_\Disc \right)
 			=\Bigl\llbracket \underline{C}_\Disc +\Pd_\Disc^{n+1}-\Nd_\Disc^{n+1}, \v_\Disc \Bigl\rrbracket_\M 
 			\quad \forall \v_\Disc \in \V_{\Disc,0}^{\Gamma_D}.\label{sch:HFV:phi}
\end{gather}
\end{subequations}
As previously, we use the notation $ \log(\u_\Disc) = \left( \left(\log(\uk)\right)_{K\in{\M}},\left( \log(u_\s)\right)_{\s \in \E}\right)$.
\vspace{-0.3cm}
\subsection{Some structural differences between schemes} 
As highlighted by the unified presentation above, both schemes are very similar and rely on the same features. 
Note that both local reconstruction operators $r^\Ds$ and $r^K$ take into account all the local unknowns of the geometric entity considered (diamond or cells), this property is the key point of the analysis of this kind of schemes, see \cite{CCHK18,Moa:22}. 
\newline
However, the schemes exhibit differences, some of which are listed below:
\vspace{-0.15cm}
\begin{itemize}	
	\item the discrete HFV gradient $\nabla_\Disc$ includes a stabilisation term for the sake of coercivity and the stabilisation parameter $\eta$ has to be chosen a priori, whereas the DDFV one is simpler and do not need any choice of parameter; 
	\item the DDFV unknowns are all "volumic", in the sense that there are associated to geometric entities with non-zero two-dimensional measures, whereas the faces unknowns of the HFV method have no mass and have no influence on the discrete time derivative terms $\left \llbracket\Nd_\Disc^{n+1}-\Nd_\Disc^{n}, \v_\Disc \right \rrbracket_\M$ and $\left \llbracket\Pd_\Disc^{n+1}-\Pd_\Disc^{n}, \v_\Disc \right \rrbracket_\M$; 
	\item the cells unknowns of the HFV scheme can be eliminated before solving linear systems, using a static condensation procedure  (see~\cite[Section 5.1.2.]{Moa:22}), whereas one has to solve a system including all primal and dual unknowns for DDFV; 
	\item the HFV scheme can be used in 3D without any modification (the edges become faces), whereas using a DDFV method in 3D requires more sophisticated changes (see \cite{CH_11}). 
\end{itemize}

\section{Numerical experiments}\label{sec:num}
%
The two numerical schemes described here are nonlinear, hence their algebraic realisations boil down to the resolution of nonlinear systems of equations.
To do so, we use Newton method, with an adaptative time stepping strategies: if the Newton method does not converge, we try to compute the solution for a smaller time step $0.5 \times \Delta t$. If the method converges, we use a bigger time step $1.4 \times \Delta t$. 
\JM{The initial time step is denoted by $\Delta t_{ini}$, and we  also impose a maximal time step $\Delta t_{max}$.} 
For the HFV scheme, at each system resolution, a static condensation is used to eliminate the cell unknowns (see~\cite[Section 5.1.2.]{Moa:22}), and we use $\eta = 1.5$.
\JM{Note that we use $N$, $P$ and $\phi$ as discrete unknowns in the schemes.}

\begin{wrapfigure}{r}{4.3cm}
	\vspace{-1.1cm}
  \centering
  	{\def\svgwidth{4.1cm}
	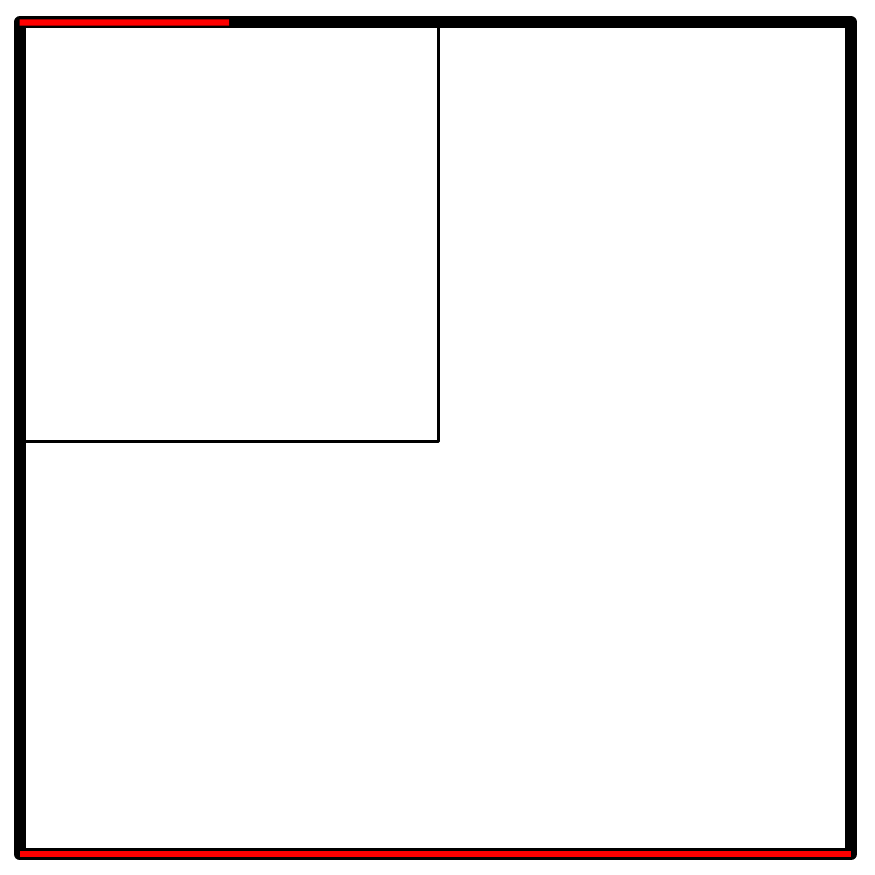
	}
	\centering
\caption{PN diode geometry.}	
\vspace{-0.3cm}
\label{fig:diode}
\end{wrapfigure}

%
The test case used below follow the framework used in \cite{Moa:22} to describe a 2D PN-junction, whose geometry is described in Figure~\ref{fig:diode}.\newline
%
%
The domain $\Omega$ is the unit square ${]0,1[}^2$. For the boundary conditions, we split $\Gamma^D = \Gamma^D_0 \cup \Gamma^D_1 $ with $\Gamma^D_0 = [0,1] \times \{0\}$ and $\Gamma^D_1 = [0, 0.25] \times \{1\}$. 
For $i \in \{0,1\}$, we let
$N^D = N^D_i$, $P^D = P^D_i$ and $\phi^D = \frac{\log(N^D_i) - \log(P^D_i)}{2}$ on $\Gamma^D_i$.
To be consistent with the compatibility condition~\eqref{eq:compcond} we assume that there exists a constant $\alpha_0$ such that 
$\displaystyle 	\log(N^D \times P^D) = \alpha_0 $.
Therefore for given $N^D$ and $\alpha_0$ we set $P^D = \frac{\e^{\alpha_0}}{ N^D}$ on $\Gamma^D$. 
Thus, one has $\alpha_N = \alpha_P = \frac{\alpha_0}{2}$.
The doping profile $C$ is piecewise constant, equal to $-1$ in the P-region and $1$ in the N-region (see Figure~\ref{fig:diode}).
Last, we use the following smooth initial conditions: 
$N_0(x,y) = N^D_1 + (N^D_0 - N^D_1) (1 - \sqrt{y})$ and $P_0(x,y) = P^D_1 + (P^D_0 - P^D_1) (1 - \sqrt{y})$.
%

%

\subsection{Positivity}
\JM{
In this section, we compare the positivity preservation of the schemes.
The test uses the following values: $\lambda = 0.05$, $N_0^D = 0.1$, $N_1^D = 1$ and $\alpha_0 = -4$. 
We perform a test on a distorted quadrangle mesh (mesh\_quad\_6 of the FVCA 8 Benchmark), with $\Delta t_{ini} = 1.4 \, 10^{-3}$ and $\Delta t_{max} = 0.1$.
We show in Figure~\ref{fig:bounds} the evolution of the minimal values of $P$ and $N$, along with the time step and the number of Newton's iterations needed to compute the solutions at a given time for each time step. The minimal values are taken on every unknowns (primal and dual cells for the DDFV scheme, cells and faces for the HFV one).
%
\begin{figure}[h]
\pgfplotsset{width=0.95\linewidth,height=0.39\linewidth,compat=newest}
\begin{minipage}[c]{1\linewidth}
\raggedright
      \begin{tikzpicture}[scale= 0.99]
        \begin{loglogaxis}[
        		xmin=1e-3,
        		xmax=1,        		
        		ymin= 3.0e-4,
        		max space between ticks=25,
            	legend style = { 
             	at={(0.975,0.94)},
              	anchor = north east,
              	tick label style={font=\tiny},
              	legend columns=2
            			},
            	ylabel=\scriptsize{Minimal values \vphantom{Time step Distance from $a$}},
            	xlabel=\scriptsize{Time},
            	xticklabel pos=top,
          	]
          \addplot[mark=pentagon*, draw=none, RubineRed] table[x=Temps,y=min_N] {tps_full};
          \addplot[mark=pentagon*, draw=none, cyan]	table[x=Temps,y=min_P] {tps_full};          	
          \addplot[mark=triangle, draw=none, RedViolet] table[x=Temps,y=min_N] {tps_full_DDFV};
          \addplot[mark=triangle, draw=none, RoyalPurple]	table[x=Temps,y=min_P] {tps_full_DDFV};
          \legend{\tiny $\min \Nd_\Disc^n$ (HFV),\tiny $\min \Pd_\Disc^n$ (HFV), \tiny $\min \Nt^n$ (DDFV),\tiny $\min \Pt^n$ (DDFV) } 
	      \end{loglogaxis}
      \end{tikzpicture} 
      \begin{tikzpicture}[scale= 0.99]  
		\pgfplotsset{
    			xmin=1e-3,
        		xmax=1,
        		legend style = { 
              	at={(0.02,0.98)},
              	anchor = north west,
              	tick label style={font=\tiny},
              	legend columns=4
            			}
			}
      	\begin{loglogaxis}[
        		ymax=0.6,
        		ymin=1e-3,
        		axis y line*=left,
			axis x line=none,
            	ylabel=\scriptsize{Time step \vphantom{Distance from $a$ Minimal values}},
          	]
          	\addplot[mark=pentagon*, mark size=1.9,  line width=0.5pt, RubineRed] table[x=Temps,y=time_step] {tps_full};\label{step}          	
          	\addplot[mark=triangle, mark size=1.9,  line width=0.5pt, draw= RedViolet] table[x=Temps,y=time_step] {tps_full_DDFV};\label{step_DDFV}
	   	\end{loglogaxis}
        	\begin{semilogxaxis}[
        		ymin = 0.,
        		axis y line*=right,
        		ymax=7,
        		max space between ticks= 20,
            	ylabel=\scriptsize{\#  of Newton iterations},
            	xlabel=\scriptsize{Time}
          	]
          	\addplot[mark=pentagon*,mark size=2.5, draw=none, cyan] table[x=Temps,y=Nb_iter] {tps_full};          
          	\addplot[mark=triangle,mark size=2.5, draw=none, RoyalPurple] table[x=Temps,y=Nb_iter] {tps_full_DDFV};
          	\addlegendimage{/pgfplots/refstyle=step};   
          	\addlegendentry{\tiny \# iterations (HFV)}; 
          	\addlegendentry{\tiny \# iterations (DDFV)};    
          	\addlegendimage{/pgfplots/refstyle=step_DDFV};    	
			\addlegendentry{\tiny Time step (HFV)} 	       	
          	\addlegendentry{\tiny Time step (DDFV)}        	
	      	\end{semilogxaxis}
      \end{tikzpicture} 
\end{minipage}
	\caption{\textbf{Test-case}. Evolution of the discrete minimal values, time step and cost}
	\label{fig:bounds}
\end{figure}
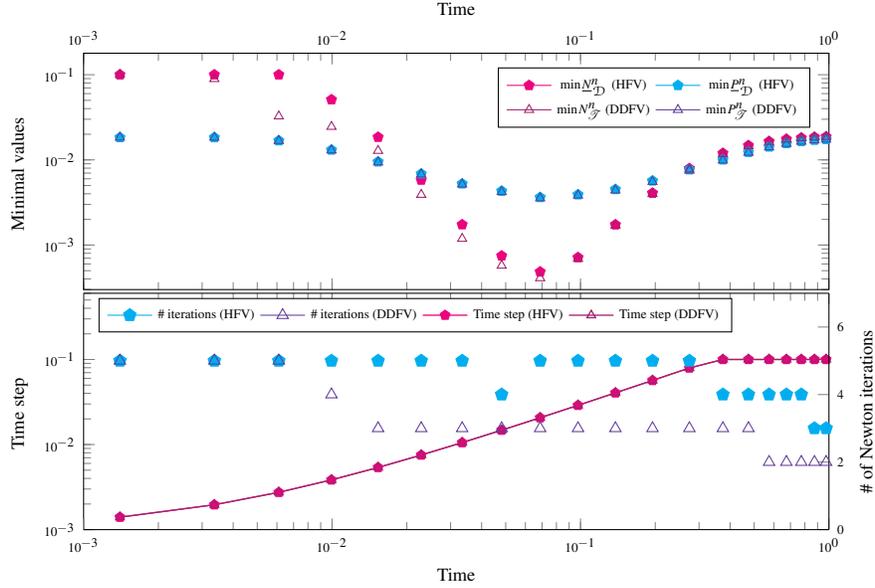 
One can see that both schemes compute, as expected by the theoretical results, positive densities. The minimal values computed are of the same order for both schemes.
Moreover, both computations proceed without the need of a time step reduction. 
Regarding the cost, it appears that the HFV scheme needs more Newton iterations than the DDFV one (90 vs 63). For both schemes, the number of iteration decay as the time increases, since the solutions converge exponentially fast towards the equilibrium.
}
\subsection{Long-time behaviour}
Here, we investigate the long-time behaviour of the schemes. 
%
At the continuous level, one usually quantify the distance between the solution $(N,P,\phi)$ and the equilibrium $(N^e,P^e,\phi^e)$ by looking at the relative entropy, defined as 
\vspace{-0.2cm}
\[
	\mathbb{E}(t) = \int_\Omega N^e H \left ( \frac{N}{N^e} \right )+ \int_\Omega P^e H \left ( \frac{P}{P^e} \right )  \\
		+ \frac{\lambda^2}{2} \| \nabla(\phi-\phi^e) \|^2_{L^2(\Omega)},
\]
with $H : s\mapsto s \log(s) -s +1$.
One can check that $(N,P,\phi)$ coincides with the equilibrium if and only if the relative entropy cancels.
In the following, we are interested in the evolution of the discrete counterparts of this quantities, 
\JM{
defined as 
\vspace*{-0.15cm}
\[
	\mathbb{E}^n_\Disc =  
	\left \llbracket \Nd_\Disc^e H \left ( \frac{\Nd_\Disc^n}{\Nd_\Disc^e}  \right), \one_\Disc \right \rrbracket_\M
	+\left \llbracket \Pd_\Disc^e H \left ( \frac{\Pd_\Disc^n}{\Pd_\Disc^e}  \right), \one_\Disc \right \rrbracket_\M
	+ \frac{\lambda^2}{2} a_\Disc \left( \phid_\Disc^{n}-\phid_\Disc^e  , \phid_\Disc^{n}-\phid_\Disc^e \right)
\]
for the HFV scheme (where $\one_\Disc$ is the discrete elements whose coordinates are $1$, and the product, quotient and functions are applied coordinate-wise) and similar definition for the DDFV scheme. 
Note that the HFV entropy does not take into account the edge unknowns of the discrete densities. 
To compute the discrete equilibrium, we use a nonlinear scheme for \eqref{pb:PoissonBoltz} and get $\phid_\Disc^e$, then we defined the associated densities following the continuous relations $N^e = e^{\alpha_N + \phi^e}$ and $P^e = e^{\alpha_P - \phi^e}$.
}
%
We consider a test case with physical data $N_0^D = \e$, $N_1^D = 1$ and $\alpha_0 = 0$.
We also use two different values of the Debye length $\lambda$, respectively $1$ and $0.01$.
\JM{We perform simulations on a triangular mesh, with a  $\Delta t_{ini} = \Delta t_{max} = 0.1$. } 
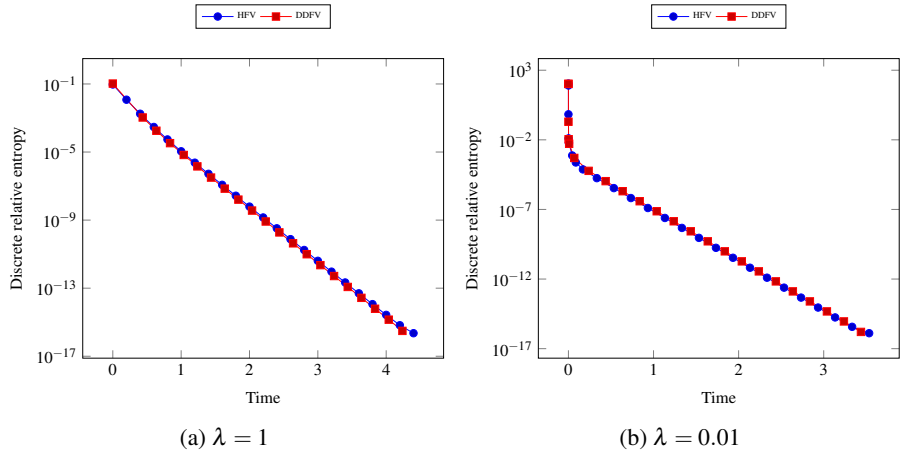
\begin{figure}[ht]
\begin{subfigure}{0.49\linewidth}
\center
\begin{tikzpicture}[scale= 0.7]
        \begin{semilogyaxis}[   
            	 legend style = { 
              at={(0.5,1.1)},
              anchor = south,
              tick label style={font=\footnotesize},
              legend columns=-1
            },ylabel=\small{Discrete relative entropy},xlabel=\small{Time}
          ]
          \addplot table[x=Temps,y=Entro] {data/HFV/mesh1_4/tps_impli_1};
          \addplot table[x=Time,y=Entro] {data/DDFV/mesh1_4/l_1};
          \legend{\tiny HFV , \tiny DDFV};
	      \end{semilogyaxis}
\end{tikzpicture}    
\caption{$\lambda = 1$}
\label{fig:tps:1}
\end{subfigure}
\begin{subfigure}{.49\linewidth}
\center
\begin{tikzpicture}[scale= 0.7]
        \begin{semilogyaxis}[
            legend style = { 
              at={(0.5,1.1)},
              anchor = south,
              tick label style={font=\footnotesize},
              legend columns=-1
            },ylabel=\small{Discrete relative entropy},xlabel=\small{Time}
          ]           
          \addplot table[x=Temps,y=Entro] {data/HFV/mesh1_4/tps_impli_0.01};
          \addplot table[x=Time,y=Entro] {data/DDFV/mesh1_4/l_0.01};
          \legend{\tiny HFV, \tiny DDFV};
        \end{semilogyaxis}
\end{tikzpicture}
\caption{ $\lambda = 0.01$ }
\label{fig:tps:0.01}
\end{subfigure}   
\caption{\textbf{Long-time behaviour.} Evolution of the discrete relative entropies.}
\label{fig:tps}
\end{figure}
On Figure \ref{fig:tps}, we show the evolutions of the discrete relative entropies along time, for the two values of $\lambda$ and both schemes.
As expected, the convergence towards the equilibrium is exponentially fast, as in the continuous framework. 
Moreover, it is remarkable to notice that the decay rates are almost the same for both schemes. 
Moreover, with the small Debye length (Figure \ref{fig:tps:0.01}), both schemes are able to capture the behaviour with a very fast evolution far from the equilibrium, then slower once close to it.
%
%
%

\end{document}

%% file: macroDDFV.tex
\newcommand{\dr}{\partial}

\def\k{{{K}}}
\def\l{{{L}}}
\newcommand{\ke}{{ {K}^*}}
\def\le{{{L}^*}}
\newcommand{\sig}{ \sigma}
\newcommand{\sige}{\sigma^*}

\newcommand{\M}{\mathfrak{M}}
\newcommand{\Mie}{\mathfrak{M}^*}

\newcommand{\ut}{u_{\T}}
\newcommand{\Pt}{P_{\T}}

\newcommand{\Nt}{N_{\T}}

\newcommand{\vt}{v_{\T}}
\newcommand{\wt}{w_{\T}}
\newcommand{\pott}{\phi_{\T}}

\newcommand{\Ee}{\mathcal{E}}
\newcommand{\Ds}{{\mathcal{D}}}

\newcommand{\nksig}{{\mathbf{n}}_{\sig K}}

\newcommand{\nkesige}{\boldsymbol{{\mathbf{n}}}_{\sige\ke}}

\newcommand{\tkele}{\boldsymbol{{\boldsymbol{\tau}}}_{\boldsymbol{\ke\!\!\!,\le}}}
\newcommand{\tkl}{\boldsymbol{{\boldsymbol{\tau}}}_{\boldsymbol{K\!,L}}}

\newcommand{\petitt}{{\scriptscriptstyle\mathcal{T}}}
\newcommand{\DD}{\mathfrak{D}}

\newcommand{\Dsig}{{\mathcal{D}_{\sigma,\sigma^*}}}

\newcommand{\uk}{u_{\k}}
\newcommand{\ul}{u_{\l}}
\newcommand{\uke}{u_{\ke}}
\newcommand{\ule}{u_{\le}}

\newcommand{\vk}{v_\k}

\newcommand{\vke}{v_\ke}

\newcommand{\xk}{x_{\k}}
\newcommand{\xl}{x_{\l}}

\newcommand{\xke}{x_{\ke}}
\newcommand{\xle}{x_{\le}}

\newcommand{\msig}{{\rm m}_{\sig}}

\newcommand{\msige}{{\rm m}_{\sige}}

\newcommand{\mk}{{\rm m}_{\k}}
\newcommand{\mke}{{\rm m}_{\ke}}
\newcommand{\md}{{\rm m}_{\Ds}}

\newcommand{\dsp}{\displaystyle}

\newcommand{\grad}{\nabla}

\newcommand{\R}{\mathbb R}

\newcommand{\gradDD}{\nabla^{\DD }}

\newcommand{\gradD}{\nabla^{\Ds}}

\newcommand{\T}{\mathcal{T}}
\newcommand{\Rt}{\R^{\petitt}}

\newcommand{\RdDD}{\left(\R^2\right)^{\DD}}

\newcommand{\sumdiam}{\ssum_{\Ds\in\DD}}

\newcommand{\sumpri}{\ssum_{\k\in\M}}
\newcommand{\sumdua}{\ssum_{\ke\in{\overline\Mie}}}

\DeclareMathOperator*{\ssum}{\sum}

\newcommand{\xib}{{\boldsymbol \xi}}
\newcommand{\phib}{{\boldsymbol \varphi}}



%% file: SC_typique.pdf_tex
\begingroup%
  \makeatletter%
  \providecommand\color[2][]{%
    \errmessage{(Inkscape) Color is used for the text in Inkscape, but the package 'color.sty' is not loaded}%
    \renewcommand\color[2][]{}%
  }%
  \providecommand\transparent[1]{%
    \errmessage{(Inkscape) Transparency is used (non-zero) for the text in Inkscape, but the package 'transparent.sty' is not loaded}%
    \renewcommand\transparent[1]{}%
  }%
  \providecommand\rotatebox[2]{#2}%
  \newcommand*\fsize{\dimexpr\f@size pt\relax}%
  \newcommand*\lineheight[1]{\fontsize{\fsize}{#1\fsize}\selectfont}%
  \ifx\svgwidth\undefined%
    \setlength{\unitlength}{419.52755906bp}%
    \ifx\svgscale\undefined%
      \relax%
    \else%
      \setlength{\unitlength}{\unitlength * \real{\svgscale}}%
    \fi%
  \else%
    \setlength{\unitlength}{\svgwidth}%
  \fi%
  \global\let\svgwidth\undefined%
  \global\let\svgscale\undefined%
  \makeatother%
  \begin{picture}(1,1)%
    \lineheight{1}%
    \setlength\tabcolsep{0pt}%
    \put(0.52925868,0.27331711){\color[rgb]{0,0,0}\makebox(0,0)[lt]{\lineheight{1.25}\smash{\begin{tabular}[t]{l}N-region \\ $C = 1$\end{tabular}}}}%
    \put(0.13280309,0.69886021){\color[rgb]{0,0,0}\makebox(0,0)[lt]{\lineheight{1.25}\smash{\begin{tabular}[t]{l}P-region \\ $C = -1$\end{tabular}}}}%
    \put(0.14142531,0.0714536){\color[rgb]{0,0,0}\makebox(0,0)[lt]{\lineheight{1.25}\smash{\begin{tabular}[t]{l}$\color{red}{\Gamma^D_0}$\end{tabular}}}}%
    \put(0,0){\includegraphics[width=\unitlength,page=1]{SC_typique.pdf}}%
    \put(0.05636858,0.88043292){\color[rgb]{0,0,0}\makebox(0,0)[lt]{\lineheight{1.25}\smash{\begin{tabular}[t]{l}$\color{red}{\Gamma^D_1}$\end{tabular}}}}%
  \end{picture}%
\endgroup%